\documentclass{amsart}
\usepackage{amsfonts,amssymb,amsmath,amsthm}
\usepackage{url}
\usepackage{enumerate}
\usepackage{euscript}
\usepackage{verbatim}

\newtheorem{thm}{Theorem}[section]
\newtheorem{lemma}[thm]{Lemma}
\newtheorem{cor}[thm]{Corollary}
\newtheorem{prop}[thm]{Proposition}

\theoremstyle{definition}

\numberwithin{equation}{section}

\newcommand{\nin}{\noindent}
\newcommand{\bn}{\bigskip\noindent}
\newcommand{\mn}{\medskip\noindent}

\newcommand{\DEF}[1]{\textit{#1\/}}
\newcommand{\raf}[1]{\/{\rm (\ref{#1})}}
\newcommand{\mystack}[2]{\stackrel{\scriptstyle{#1}}{#2}}

\newcommand{\C}[2]{{{#1}\choose{{#2}}}}

\newcommand{\emp}{\varnothing}
\newcommand{\ddef}{\mbox{$\: :=\:$}}
\newcommand{\deff}{\mbox{$\: =:\:$}}
\newcommand{\sm}{\smallsetminus}
\newcommand{\sub}{\subseteq}
\newcommand{\ZZ}{\mathbb{Z}}
\newcommand{\QQ}{\mathbb{Q}}
\newcommand{\RR}{\mathbb{R}}

\newcommand{\beq}[1]{\begin{equation}\label{#1}}
\newcommand{\enq}{\end{equation}}

\newcommand{\eps}{\varepsilon}

\newcommand{\DD}{\mathcal{D}}

\newcommand{\EB}{\EuScript{B}}

\newcommand{\Di}{{\mathcal{D}_{1}}}
\newcommand{\Dii}{{\mathcal{D}_{2}}}
\newcommand{\Diii}{{\mathcal{D}_{3}}}
\newcommand{\Diiii}{{\mathcal{D}_{4}}}
\newcommand{\DC}{{\mathcal{D}'}}
\newcommand{\nge}{{\lceil n^{g\eps}\rceil}}
\newcommand\chic{\chi_{\rm c}}

\author{Ararat Harutyunyan \and P. Mark Kayll \and Bojan Mohar \and Liam Rafferty}

\address{
Department of Mathematics \\
Simon Fraser University \\
Burnaby, B.C. V5A~1S6, Canada}
\email{aha43@sfu.ca}
\thanks{AH: This work forms part of the author's PhD 
dissertation~\cite{Harutyunyan2011}; research supported by
FQRNT (Le Fonds qu\'{e}b\'{e}cois de la recherche sur la nature et
les technologies) doctoral scholarship.}

\address{
Department of Mathematical Sciences \\
University of Montana \\
Missoula MT 59812-0864, USA}
\email{mark.kayll@umontana.edu}
\thanks{PMK: Contact author; research initiated, in part,
under UM's sabbatical program.}

\address{
Department of Mathematics \\
Simon Fraser University \\
Burnaby, B.C. V5A~1S6, Canada}
\email{mohar@sfu.ca}
\thanks{BM: Supported in part by an NSERC Discovery Grant
(Canada),  by the Canada Research Chair program, and by the
  Research Grant P1--0297 of ARRS (Slovenia). On leave from:
  IMFM \& FMF, Department of Mathematics, University of Ljubljana, Ljubljana,
  Slovenia.}

\address{
Department of Mathematical Sciences \\
University of Montana \\
Missoula MT 59812-0864, USA}
\email{rafferty@member.ams.org}
\thanks{LR: This work forms part of the author's PhD 
dissertation~\cite{Rafferty2011}; research supported in part by a UM
Graduate Student Summer Research Award funded by the George and
Dorothy Bryan Endowment.\\[0.5em]
Copyright \copyright\ 2011 by A. Harutyunyan, P.M. Kayll, B. Mohar and L. Rafferty}

\begin{document}

\setlength{\baselineskip}{1.247 \baselineskip}
\enlargethispage*{1em}

\title[Uniquely {\boldmath$D$}-colourable digraphs with large girth]{Uniquely {\boldmath$D$}-colourable digraphs with large girth}

\begin{abstract}
  Let $C$ and $D$ be digraphs. A mapping $f:V(D)\to V(C)$ is a
  $C$-colouring if for every arc $uv$ of $D$, either $f(u)f(v)$
  is an arc of $C$ or $f(u)=f(v)$, and the preimage of every
  vertex of $C$ induces an acyclic subdigraph in $D$. We say
  that $D$ is $C$-colourable if it admits a $C$-colouring and
  that $D$ is uniquely $C$-colourable if it is surjectively
  $C$-colourable and any two $C$-colourings of $D$ differ by an
  automorphism of $C$. We prove that if a digraph $D$ is not
  $C$-colourable, then there exist digraphs of arbitrarily large
  girth that are $D$-colourable but not
  $C$-colourable. Moreover, for every digraph $D$ that is
  uniquely $D$-colourable, there exists a uniquely
  $D$-colourable digraph of arbitrarily large girth. In
  particular, this implies that for every rational number $r\geq
  1$, there are uniquely circularly $r$-colourable digraphs with
  arbitrarily large girth.
\end{abstract}
\maketitle

\section{Introduction}

In a seminal \emph{Canadian Journal of Mathematics}
article~\cite{Erdos59}, Paul Erd\H{o}s established
nonconstructively the existence of graphs with arbitrarily large
girth $\gamma$ and arbitrarily large chromatic number $\chi$.
In these introductory remarks, we focus mainly on 
Erd\H{o}s' theorem, for the features that make it interesting are
shared by its progeny 
(e.g.\ \cite{Bollobas-Sauer76,Zhu-Nesetril,Zhu-uniq}), the first
of which also appeared in \emph{CJM}.

Because graphs with similar combined properties as guaranteed by
Erd\H{o}s  had been constructed earlier---e.g.\
triangle-free plus large $\chi$ 
\cite{Descartes47,Mycielski55,Zykov49} or
girth at least six plus large $\chi$ 
\cite{Descartes54,KellyKelly54}---his result was not wholly
unanticipated. Nevertheless, it remains somewhat
counterintuitive. Na\"{\i}vely, one might reason that having
large girth implies edge-sparsity, while having large chromatic
number entails edge-abundance, so how could both properties
coexist in one graph? Even upon closer inspection, the result
seems paradoxical. If, for a positive integer $\ell$, a graph
$G$ satisfies $\gamma > \ell$, then any set of at most $\ell$
vertices induces an acyclic, hence $2$-colourable, subgraph. Why
is it not possible to assemble such colourings into a proper
colouring of $V(G)$ using few colours? These questions' answers
might be regarded as the take-home message of Erd\H{o}s' theorem:
since the chromatic number depends intrinsically upon the
graph's global structure, local $2$-colourability imparts
nothing on $\chi$.

Yet the theorem's influence somehow manages to transcend its
important message. Every student of combinatorial probability
studies Erd\H{o}s' proof (cf.\
\cite{AignerZiegler2010,Alon-Spencer2000,Bollobas98,BondyMurty08,Diestel2010,MolloyReed2002}), 
that employs a
deterministic step after a probabilistic argument for the
existence of graphs with \emph{few} short cycles and small
stability number. And though constructive
proofs~\cite{Kriz,Lovasz68,NesetrilRodl79} of Erd\H{o}s'
existence theorem eventually followed, their complexity perhaps
precludes their inclusion in `The Book', generally imagined to
favour the elegance, clarity, and simplicity of Erd\H{o}s'
original argument.

Aside from its beautiful proof, the theorem's influence can also
be measured by considering its descendants.
Ne\v{s}et\v{r}il~\cite{Nesetril73} conjectured, and Bollob\'{a}s
and Sauer~\cite{Bollobas-Sauer76} proved, the existence of
graphs as guaranteed by Erd\H{o}s that are, moreover, uniquely
$\chi$-colourable. Colourings are special cases of homomorphisms
into a fixed graph, and Zhu~\cite{Zhu-uniq} extended both
Erd\H{o}s' and Bollob\'{a}s and Sauer's results to homomorphisms 
into general graphs. Rather recently, the results of \cite{Zhu-uniq} were
extended by Ne\v{s}et\v{r}il and Zhu~\cite{Zhu-Nesetril} to
give a simultaneous generalization of Zhu's two primary
results. Without attempting to give an exhaustive list, we also
note the appearance in recent years of a host of other articles
related to the interplay between girth and colouring; see, e.g.,
\cite{EHK98,KostochkaNesetril99,LinZhu2006,Muller79,Nesetril-Zhu2001,PanZhu2005,Zhu99}.
The results of the present paper extend the main theorems of
Zhu~\cite{Zhu-uniq} to digraphs with acyclic homomorphisms.

\subsection*{Notation, terminology, details}

As much as possible, we try to follow standard terminology. See,
for example, \cite{BG2009,BondyMurty08} for graphs and digraphs,
\cite{Alon-Spencer2000,MolloyReed2002} for probabilistic
concerns, and \cite{HellNesetril2004} for homomorphisms.

Our digraphs are simple---i.e.\ loopless and without
multiple arcs---however, we allow two vertices $u$, $v$ to be joined
by two oppositely directed arcs, $uv$ and $vu$. The \DEF{girth}
of a graph or digraph refers to the length of a shortest cycle,
that we take to mean directed cycle in the digraph case (and
infinite in either acyclic case).

Recall that a \DEF{homomorphism} of a graph $G$ into a graph $H$
is a function $\phi\colon V(G) \to V(H)$ such that
$\{\phi(u),\phi(v)\}\in E(H)$ whenever $\{u,v\}\in E(G)$. 
An \DEF{acyclic homomorphism} of a digraph $D$ into a digraph
$C$ is a function $\phi\colon V(D) \to V(C)$ such that:
\begin{itemize}
\item[(i)] for every vertex $v\in V(C)$, the subdigraph of $D$ induced
by $\phi^{-1}(v)$ is acyclic;
\item[(ii)] for every arc $uv\in E(D)$, either $\phi(u)=\phi(v)$, or
$\phi(u)\phi(v)$ is an arc of $C$.
\end{itemize}
If digraphs $D$ and $C$ are obtained from undirected graphs $G$ and $H$,
respectively, by replacing every edge by two oppositely directed arcs,
then acyclic homomorphisms between $D$ and $C$ correspond to usual graph
homomorphisms between $G$ and $H$. In this sense, acyclic homomorphisms
can be viewed as a generalization of the notion of homomorphisms of
undirected graphs.

It is well-known and easy to see that a graph $G$ is (properly)
$r$-colourable (for a positive integer $r$) if and only if $G$
admits a homomorphism to the complete graph $K_r$. Thus, $G$ is
commonly called \DEF{$H$-colourable} if there is a
homomorphism from $G$ to $H$. In the same way as homomorphisms
generalize the notion of graph colouring,  acyclic
homomorphisms generalize digraph colouring;
cf.~\cite{BFJKM2002}. Motivated by this, we say that 
a digraph $D$ is \DEF{$C$-colourable} if there is an acyclic
homomorphism from $D$ to $C$.

Zhu generalized Erd\H{o}s' theorem as follows.
\begin{thm}[\cite{Zhu-uniq}]
\label{Zhu-one}
If $G$ and $H$ are graphs such that $G$ is not $H$-colourable,
then for every positive integer $g$, there exists a graph $G^*$ of girth
at least $g$ that is $G$-colourable but not $H$-colourable.
\end{thm}
\noindent
To recover Erd\H{o}s' theorem, suppose that we want to arrange
for $\gamma\geq g$ and $\chi\geq r$, for some prescribed integers
$g$ and $r$; then we take $G=K_r$ and $H=K_{r-1}$ in
Theorem~\ref{Zhu-one}. 

Our first main result is a digraph analogue of the preceding
result.
\begin{thm}
\label{main-one}
If $D$ and $C$ are digraphs such that $D$ is not $C$-colourable,
then for every positive integer $g$, there exists a digraph $D^*$ of girth
at least $g$ that is $D$-colourable but not $C$-colourable.
\end{thm}
\noindent
Just as Theorem~\ref{Zhu-one} generalizes Erd\H{o}s' theorem,
Theorem~\ref{main-one} generalizes the analogue appearing in 
\cite{BFJKM2002}. See the introduction to
Section~\ref{sect:chic} for a statement of this analogue.

A graph $G$ is \DEF{uniquely $H$-colourable} if it is
surjectively $H$-colourable, and for any two $H$-colourings 
$\phi$, $\psi$ of
$G$, there is an automorphism $\pi$ of $H$ such that 
\begin{equation}
\label{differ-aut}
\phi = \pi\circ\psi.
\end{equation}
Unique $D$-colourability is defined analogously for digraphs
$D$. In either case, when \raf{differ-aut} occurs, we sometimes
say that $\phi$ and $\psi$ \DEF{differ by an automorphism} of $H$. 
A graph $H$ is a \DEF{core} if it is uniquely $H$-colourable;
likewise for digraphs. To align this formulation with the usual
one (cf.~\cite{GodsilRoyle2001,HellNesetril2004}), we offer the
following observation about the digraph version.

\begin{lemma}
\label{lem:core}
A digraph $D$ is a core if and only if every acyclic homomorphism
$V(D)\to V(D)$ is a bijection.
\end{lemma}

\proof
Let $\phi:V(D)\to V(D)$ be an acyclic homomorphism. If $\phi$ is not
a bijection, then $\phi$ and the identity homomorphism do not differ by an
automorphism of $D$, so $D$ is not a core.

Suppose now that $D$ is not a core, and let $\phi,\psi$ be two acyclic
homomorphisms that do not differ by an automorphism of $D$.
If $\phi$ (or $\psi$) is bijective, then it is a homomorphism of $D$
onto itself. This implies that it is an automorphism of $D$. Therefore,
$\phi$ and $\psi$ are not both bijective.
\qed

Zhu generalized the aforementioned Bollob\'{a}s-Sauer
theorem~\cite{Bollobas-Sauer76} as follows. 
\begin{thm}[\cite{Zhu-uniq}]
\label{Zhu-two}
For every graph $H$ that is a core and every positive integer
$g$, there exists a graph $H^*$ of girth at least $g$ that is
uniquely $H$-colourable. 
\end{thm}
\noindent
Bollob\'{a}s and Sauer's result follows from
Theorem~\ref{Zhu-two} because complete graphs are cores, as is
easily verified.

Our second main result establishes a digraph analogue of
Theorem~\ref{Zhu-two}. 
\begin{thm}
\label{main-two}
For every digraph $D$ that is a core and every positive integer
$g$, there exists a digraph $D^*$ of girth at least $g$ that is
uniquely $D$-colourable.
\end{thm}
\noindent
Theorem~\ref{main-two} immediately applies to digraph colourings
and digraph circular colourings (see, e.g., \cite{BFJKM2002,Mo})
to yield our third main result, Theorem~\ref{thm:colouring}. In
favour of an abbreviated mention of this result here, we
postpone until Section~\ref{sect:chic} its full statement, the
definition of `circular colouring', and some related
discussion.  
\begin{cor}
\label{cor:colouring}
For every rational number $r\geq 1$ and every positive integer
$g$, there exists a digraph of girth at least $g$ that is
uniquely circularly $r$-colourable.
\end{cor}

We devote Section~\ref{sect:proof2}
to the proof of Theorem~\ref{main-one}, while
Section~\ref{sect:proof3} contains the proof of
Theorem~\ref{main-two}.  Both proofs are probabilistic and follow
the main ideas of \cite{Bollobas-Sauer76} and \cite{Zhu-uniq},
which themselves trace back to Erd\H{o}s~\cite{Erdos59}.
However,  just as both of these earlier refinements required new
ideas to move to the next level, additional care and some
inspiration are needed to extend the 
proofs to the digraph setting.

\section{Proof of Theorem~\protect{\ref{main-one}}}

\label{sect:proof2}

We begin by setting up a suitable random digraph model.
Suppose that $V(D) = \{1,2,\dots,k\}$ and that
$q=|E(D)|$. Let $n$ be a (large) positive integer, and let $D_n$
be the digraph obtained from $D$ as follows: replace every vertex
$i$ with a (temporarily) stable set $V_i$ of $n$ ordered vertices
$v_1,v_2,...,v_n$, and replace each arc $ij$ of $D$ by the set of
all possible $n^2$ arcs from $V_i$ to $V_j$; additionally, add
each arc $v_rv_s$ such that $v_r,v_s\in V_i$ and $r<s$. Clearly,
$|V(D_n)| = kn$ and $|E(D_n)| = qn^2+k\C{n}{2}$.

Now fix a positive $\eps < 1/(4g)$. Our random digraph model
$\DD=\DD(D_n,p)$ consists of those spanning subdigraphs of $D_n$
in which the arcs of $D_n$ are chosen randomly and independently
with probability $p=n^{\eps-1}$.

As usual in nonconstructive probabilistic proofs of results of
this nature (cf.\ \cite{Bollobas-Sauer76,Zhu-Nesetril,Zhu-uniq}),
the idea is to show that most digraphs in $\DD$ have
only a few short cycles, and for most digraphs $H\in\DD$, the
subdigraph of $H$ obtained by removing an arbitrary yet small set of
arcs is not $C$-colourable. Choosing an $H\in\DD$ with both these
properties, we can force the girth to be large by deleting an arc
from each short cycle. Since the set $A_0$ of deleted arcs is
small, the resulting digraph $H - A_0$ satisfies the desired
conclusion of Theorem \ref{main-one}.

To make this description more precise, let $\Di$ denote the set
of digraphs in $\DD$ containing at most $\nge$ cycles of length
less than $g$, and let $\Dii$ be the set of digraphs $H\in\DD$
that have the property that $H - A_0$ is not $C$-colourable for
any set $A_0$ of at most $\nge$ arcs. We will show that
\beq{D1-lb}
\left|\Di\right|>\left(1-n^{-\eps/2}\right)\left|\DD\right|
\enq
and
\beq{D2-lb}
\left|\Dii\right|>\left(1-e^{-n}\right)\left|\DD\right|.
\enq
Since \raf{D1-lb} and \raf{D2-lb} imply that
$\Di\cap\Dii\neq\emp$ (for sufficiently large $n$), there exists
a digraph $H\in\Di\cap\Dii$. Now $H\in\Di$ implies that there is
a set $A_0$ of at most $\nge$ arcs whose removal leaves a
digraph $D^*\ddef H - A_0$ of girth at least $g$, while
$H\in\Dii$ means that $D^*$ is not $C$-colourable. Thus, it
remains to establish \raf{D1-lb} and \raf{D2-lb}.

\bn
\textit{Proof of\/} \raf{D1-lb}. The expected number $N_{\ell}$
of cycles of length $\ell$ in a digraph $H\in\DD$ is at most
\begin{equation}
\label{eq:Nl}
\C{kn}{\ell}(\ell-1)!\,p^{\ell}
\end{equation}
since there are $\C{kn}{\ell}(\ell-1)!$ ways of choosing a
cyclic sequence of $\ell$ vertices as a candidate for a cycle,
and such an  $\ell$-cycle occurs in $\DD$ with probability
either $0$ or $p^{\ell}$. It is easy to see that the product of
the first two factors in (\ref{eq:Nl}) is smaller than
$(kn)^{\ell}/\ell$.
Therefore, if $n$ is large enough, then
\[
   \sum_{\ell=2}^{g-1}N_{\ell} \leq
   \sum_{\ell=2}^{g-1} \frac{(kn^\eps)^{\ell}}{\ell} <
   k^{g-1}n^{(g-1)\eps} < n^{-\eps/2}n^{g\eps}.
\]
Now \raf{D1-lb} follows easily from Markov's Inequality.
\qed

\bn \textit{Proof of\/} \raf{D2-lb}. We shall argue that
$|\DD\sm\Dii| < e^{-n}|\DD|$. If $H\in\DD\sm\Dii$, then there is a
set $A_0$ of at most $\nge$ arcs of $H$ so that $H - A_0$ admits
an acyclic homomorphism $h$ to $C$. Let $k'=|V(C)|$. By the
pigeonhole principle, for each $i\in V(D)$, there exists a vertex
$x_i \in V(C)$ such that $|V_i\cap h^{-1}(x_i)|\geq n/k'$. Define
$\phi\colon V(D)\to V(C)$ by setting $\phi(i)=x_i$. Since $n/k'\gg
n^{g\eps}$, the set $V_i\cap h^{-1}(x_i)$ contains a subset $W_i$
of cardinality $w:=\lceil n/(2k')\rceil$ such that no arc in $A_0$
has an end vertex in $W_i$.

Since $D$ is not $C$-colourable, the function $\phi$ is not an
acyclic homomorphism. Therefore, either there is an arc $ij\in
E(D)$ such that $\phi(i)\neq\phi(j)$ and $\phi(i)\phi(j)$ is not
an arc of $C$, or there is a vertex $v\in V(C)$ such that the
subdigraph of $D$ induced on $\phi^{-1}(v)$ contains a cycle.

We first consider the case when $ij$ is an arc of $D$ such that
$\phi(i)\neq\phi(j)$ and $\phi(i)\phi(j)$ is not an arc of
$C$. Since $h$ is an acyclic homomorphism, there are no arcs
from $W_i$ to $W_j$ in $H - A_0$. By the definition of $W_i$ and
$W_j$, neither are there such arcs in $H$.

Let us now estimate the expected number $M$ of pairs of sets
$A\subseteq V_i$, $B\subseteq V_j$, with $|A|=|B|=w$, such that
$ij\in E(D)$ and such that there is no arc from $A$ to $B$ in
$H\in\DD$ (we call such a pair $A,B$ a \emph{bad pair}). By the
linearity of expectation, we have
\begin{equation}
\label{eq:N2} M = q\C{n}{w}^2(1-p)^{w^2}<
q\Big(\frac{n^w}{w!}\Big)^2 (1-p)^{w^2}=
\frac{q(n^2(1-p)^w)^w}{(w!)^2}.
\end{equation}
Since $w$ grows no more (or less) than linearly with $n$, for
sufficiently large $n$ we have
$$n^2(1-p)^w<e^{-2k'}\mbox{ ~~and~~ } \frac{q}{(w!)^2}<\frac{1}{2}.$$
Therefore, Markov's Inequality and (\ref{eq:N2}) yield
\beq{badpair}
\Pr(\exists\mbox{ a bad pair})< \frac{e^{-n}}{2}.
\enq

Suppose now that there is a vertex $v\in V(C)$ such that $D$ contains a
cycle $Q$ whose vertices are all in $\phi^{-1}(v)$. Suppose that
$Q=i_1i_2\cdots i_t$. Observe that $2\le t\le k$. Since
$\phi(Q)=\{v\}$, we conclude that $h(W_{i_1})=h(W_{i_2})=\cdots
=h(W_{i_t})=\{v\}$. Since $h$ is an acyclic homomorphism, the
subdigraph of $H$ induced on $W_{i_1}\cup W_{i_2}\cup\cdots\cup W_{i_t}$ is
acyclic.

Let us consider all sequences of sets
$U_{j_1},U_{j_2},\ldots,U_{j_\ell}$ such that, for
$r=1,2,\ldots,\ell$, we have $U_{j_r}\subseteq V_{j_r}$ and
$|U_{j_r}|=w$, and the vertex sequence $j_1j_2\cdots j_\ell$ is a
cycle in $D$. Let $U(\ell)$ denote the subdigraph of $H$ induced on
$U_{j_1}\cup U_{j_2}\cup\cdots\cup U_{j_\ell}$, and let
$P_{\ell}:=\Pr(U(\ell)\mbox{ is acyclic})$. We call this
sequence of sets \emph{bad} if $U(\ell)$ is acyclic. Since the expected
number $N$ of bad sequences is the sum of the corresponding
expectations over all possible cycle lengths, we have
\begin{equation}
\label{expected:N}
N \leq \sum_{\ell=2}^{k} \C{k}{\ell} (\ell-1)! {\C{n}{w}}^{\ell} P_{\ell}.
\end{equation}

In order to bound $N$, we first bound the probabilities
$P_{\ell}$.
\begin{lemma}
\label{lem:P_l}
There exists a constant $\gamma>0$ (not depending on $n$) such that
$P_{\ell}\leq e^{-\gamma n^{1+\eps}}$ for every integer $\ell\in\{2,3,\ldots,k\}$.
\end{lemma}

\nin
We present two proofs of Lemma~\ref{lem:P_l}. The second invokes
the Janson Inequalities (see, e.g.,
\cite[Chapter~8]{Alon-Spencer2000}). The first uses only
elementary methods and relies in the beginning on the following
observation.  
\begin{lemma}
\label{lem:liam} A digraph $D$ is acyclic if and only if every
induced subdigraph contains a vertex of outdegree $0$.
\end{lemma}

\proof If $D$ is acyclic, then every induced subdigraph of $D$ must
be acyclic and therefore must contain a vertex of outdegree $0$.
If $D$ is not acyclic, then it must contain a cycle, the vertex
set of which induces a subdigraph containing no vertex of
outdegree $0$. \qed

\mn
\textit{Proof 1 of Lemma~\ref{lem:P_l}}.
Let $E_0$ be certain ($\Pr(E_0)=1$), and let $E_j$ be the event
that all induced subdigraphs of $U(\ell)$ with more than
$\ell w - j$ vertices have a vertex of outdegree $0$ (the outdegree in
the induced subdigraph). Lemma~\ref{lem:liam} shows that
\beq{eq:sequence}
P_{\ell}=\Pr\Big(\bigcap_{j=0}^{\ell w }
E_j\Big)=\prod_{j=0}^{\ell w - 1} \Pr(E_{j+1}|E_{j}) \leq
\prod_{j=0}^{w-1} \Pr(E_{j+1}|E_{j}).
\enq
We will call a set $S\subseteq V(U(\ell))$ an \DEF{acyclic-sink
  set} if the induced subdigraph $U(\ell)[S]$ is acyclic and
there are no arcs in $U(\ell)$ from $S$ to $V(U(\ell))\sm  S$
(so $S$ acts as a sink in $U(\ell)$).

\mn
\textbf{Claim 1:} The union of two acyclic-sink sets in
$U(\ell)$ is an acyclic-sink set in $U(\ell)$.

\mn
\textit{Proof of claim.} Let $A$ and $B$ be two acyclic-sink
sets in a digraph $U(\ell)$. Since $A$ and $B$ are both sinks in
$U(\ell)$, their union $A\cup B$ is a sink because there are no
arcs from $A\cup B$ to $V(U(\ell))\sm (A\cup B)$.
Consider the three sets $A\sm B$, $B\sm A$, and $A\cap B$; each
is a subset of an acyclic-sink set so each induces an acyclic
digraph. Since $A$ is a sink in $U(\ell)$, there can be no arcs
from $A\cap B$ to $B\sm A$. Likewise $B$ is a sink in $U(\ell)$,
so there can be no arcs from $A\cap B$ to $A\sm B$. Therefore,
$A\cup B$ induces an acyclic digraph and is consequently an
acyclic-sink set in $U(\ell)$. \qed

\mn
\textbf{Claim 2:} There exists an acyclic-sink set $S\subseteq
V(U(\ell))$ of cardinality $j$ if and only if $E_j$ occurs.

\mn
\textit{Proof of claim.} If there exists an acyclic-sink set of
cardinality $j$, then a subdigraph of $U(\ell)$ with more than
$\ell w-j$ vertices must have a nonempty intersection with
it. Any subdigraph that has nonempty intersection with an
acyclic-sink set induces a subdigraph containing a vertex of
outdegree zero.

If there is no acyclic-sink set of cardinality $j$, then the
largest acyclic-sink set is an $S'\subseteq V(U(\ell))$ such
that $|S'|<j$. Then $U(\ell)-S'$ is a subdigraph of $U(\ell)$
with cardinality greater than $\ell w-j$ and with no vertices of
outdegree $0$ (otherwise we could have added them to $S'$ and
had a larger acyclic-sink set). \qed

\mn
\textbf{Claim 3:} If $U(\ell)$ has an acyclic-sink set of
cardinality $j$, then it has an acyclic-sink set of cardinality
$j-1$.

\mn
\textit{Proof of claim.}
Suppose that $S$ is an acyclic-sink set in $U(\ell)$ of
cardinality $j$. Then the subdigraph $U(\ell)[S]$ is acyclic, so
there must be a vertex $v$ with indegree $0$ in
$U(\ell)[S]$. Consider the set $S\sm \{v\}$; this induces an
acyclic subdigraph of $U(\ell)$ because it is a subdigraph of an
acyclic digraph. There were no arcs from $S$ to $V(U(\ell))\sm
S$, and there are no arcs from $S\sm \{v\}$ to $v$, so $S\sm
\{v\}$ is a sink in $U(\ell)$. Therefore, there exists an
acyclic-sink set in $U(\ell)$ of cardinality $j-1$. \qed

\medskip

We now fix $j$ in order to estimate $\Pr(E_{j+1}|E_j)$. Let
$I=\big\{1,2,\ldots,\C{\ell w}{j}\big\}$ and let
$\big{\{}S_i\big{\}}_{i\in I}$ be the $j$-subsets of the $\ell
w$ vertices of $U(\ell)$ (in some fixed order). Let $B_i$ be the
event that $S_i$ is an acyclic-sink set in $U(\ell)$. By Claim
1, if more than one $B_i$ occurs, there must be an acyclic-sink
set of cardinality at least $j+1$, and so by Claim 3, there
exists one of cardinality exactly $j+1$. Therefore by Claim 2,
\beq{eq:intersection}
\Pr\Big(E_{j+1}|\bigcap_{i\in Y} B_i\Big) = 1 \mbox{ whenever }
Y\subseteq I \mbox{ and } |Y| \geq 2.
\enq
Now additionally fix a $B_i$, and we will estimate
$\Pr(E_{j+1}|B_i)$. Let $F$ be the event that $U(\ell)-S_i$
contains a vertex of outdegree $0$. Then
\beq{eff}
\Pr(E_{j+1}|B_i)=\Pr(E_{j+1}| F \cap B_i)\Pr(F|B_i)+\Pr(E_{j+1}|
F^C \cap B_i)\Pr(F^C|B_i).
\enq
The event $E_{j+1}$ occurs when all subsets of $V(U(\ell))$ of
cardinality greater than $\ell w - (j+1)$ induce a subdigraph in
$U(\ell)$ that has a vertex of outdegree $0$. Clearly
$U(\ell)-S_i$ has cardinality $\ell w-j$, while $F^C$ is the
event that this set induces a subdigraph with no vertex of
outdegree zero. Thus $\Pr(E_{j+1}| F^C \cap B_i) = 0$. All sets
of cardinality
exceeding $\ell w - (j+1)$ that are distinct from $V(U(\ell))\sm
S_i$ have a nonempty intersection with $S_i$, which (given
$B_i$) is an acyclic-sink set in $U(\ell)$. Therefore,
subdigraphs of $U(\ell)$ induced on these sets have a vertex of
outdegree $0$, so that $\Pr(E_{j+1}| B_i \cap F)=1$. Using these
observations, \raf{eff} reduces to
$\Pr(E_{j+1}|B_i)=\Pr(F|B_i)$. The event
$F$ is independent of the event $B_i$ since the vertices in
$S_i$ do not affect the outdegree of vertices in the subdigraph
induced by $V(U(\ell))\sm S_i$. Therefore,
$\Pr(E_{j+1}|B_i)=\Pr(F)$.

Now we estimate the probability of $F$. The probability that any
particular vertex of $U(\ell)- S_i$ has outdegree $0$ in the
induced subdigraph is bounded from above by $(1-p)^{(w - j)}$. Since
these outdegree computations are independent for each vertex,
the probability that all vertices have outdegree greater than
$0$ is bounded from below by $(1-(1-p)^{(w - j)})^{(\ell w - j)}$, so
that
\begin{eqnarray}
\Pr(E_{j+1}|B_i)=\Pr(F)&\leq& 1-((1-(1-p)^{(w - j)})^{(\ell w - j)}) \nonumber \\
\label{eq:pprime}
&<& (\ell w - j)(1-p)^{(w - j)}=:p_j.
\end{eqnarray}
\noindent
We also need to estimate $\Pr(E_{j+1}|E_j)$. By Claim 2, $E_j$
occurs if and only if $\bigcup_{i \in I}B_i$ occurs. Thus we may
rewrite $\Pr(E_{j+1}|E_j)$ using inclusion-exclusion:
\begin{eqnarray*}
\Pr(E_{j+1}|E_j)&=&\Pr\Big(E_{j+1}\big|\bigcup_{i\in I}B_i\Big)\\
&=&\frac{\Pr\Big(E_{j+1} \cap 
\big(\bigcup_{i\in I} B_i\big)\Big)}{\Pr\big(\bigcup_{i\in I}B_i\big)}\\
&=&\frac{\Pr\Big(\bigcup_{i\in I}( E_{j+1} \cap B_i)\Big)}{\Pr\big(\bigcup_{i\in
I}B_i\big)}\\
&=&\sum_{\emp\neq Y\sub I} (-1)^{|Y|+1}
\frac{\Pr\Big(E_{j+1}\cap \big(\bigcap_{y\in Y}
    B_y\big)\Big)} {\Pr\big(\bigcup_{i \in I} B_i\big)}\\
&=&\sum_{\emp\neq Y\sub I} (-1)^{|Y|+1}
\frac{\Pr\Big(E_{j+1}\cap \big(\bigcap_{y\in Y}
    B_y\big)\Big)}{\Pr\big(\bigcap_{y\in Y} B_y\big)}\frac{
    \Pr\big(\bigcap_{y\in Y} B_y\big)} {\Pr\big( \bigcup_{i \in I}
    B_i\big)}\\
&=&\sum_{\emp\neq Y\sub I} (-1)^{|Y|+1}
\Pr\Big(E_{j+1}\big|\bigcap_{y\in Y} B_y\Big)\Pr\Big(\bigcap_{y\in Y}
B_y\big|\bigcup_{i \in I} B_i\Big)\\
&=&\sum_{y \in I}\Pr(E_{j+1}|B_y)\Pr\Big(B_y\big|\bigcup_{i \in I}
B_i\Big) \\
& & ~~~~~~~~ + \sum_{\mystack{Y \sub I}{|Y|\geq 2}}
(-1)^{|Y|+1}\Pr\Big(E_{j+1}\big|\bigcap_{y\in Y} B_y\Big)\Pr\Big(\bigcap_{y\in
  Y} B_y\big|\bigcup_{i\in I} B_i\Big).
\end{eqnarray*}
Using (\ref{eq:intersection}) and (\ref{eq:pprime}) in the last
expression for $\Pr(E_{j+1}|E_j)$ gives
\begin{eqnarray*}
\Pr(E_{j+1}|E_j)&\leq& p_j \sum_{y \in I}\Pr\Big(B_y\big|\bigcup_{i \in
  I} B_i\Big)+\sum_{\mystack{Y \sub I}{|Y|\geq 2}} (-1)^{|Y|+1}
\Pr\Big(\bigcap_{y\in Y} B_y\big|\bigcup_{i\in I} B_i\Big)\\
&=&p_j\sum_{y \in I}\Pr\Big(B_y\big|\bigcup_{i\in I}
B_i\Big)+\Big[\Pr\Big(\bigcup_{i\in I} B_i\big|\bigcup_{i\in I} B_i\Big)-\sum_{y \in
  I} \Pr\Big(B_y\big|\bigcup_{i\in I} B_i\Big)\Big]\\
&=&p_j\sum_{y \in I}\Pr\Big(B_y\big|\bigcup_{i\in I} B_i\Big)+\Big[1-\sum_{y \in
  I} \Pr\Big(B_y\big|\bigcup_{i\in I} B_i\Big)\Big].
\end{eqnarray*}
Since $\sum_{y \in I}\Pr\big(B_y|\bigcup_{i\in I} B_i\big)\geq 1$ 
and $p_j-1<0$, we have
\[
\Pr(E_{j+1}|E_j)~~\leq~~ 1+\sum_{y \in I}\Pr\Big(B_y\big|\bigcup_{i\in I} B_i\Big)(p_j-1)< p_j.\\
\]
\noindent 
Applying this last estimate to (\ref{eq:sequence}) yields
\begin{eqnarray}
P_{\ell} &\leq& \prod_{j=0}^{w-1} p_j=\prod_{j=0}^{w-1}
(\ell w - j) (1-p)^{(w-j)} \nonumber \\
&<&(\ell w)^w (1-p)^{w(w+1)/2} \nonumber \\
&\leq& (\ell w)^w (1-p)^{w^2 / 2} \nonumber \\
&\leq& (\ell w)^w e^{-p w^2 / 2} \nonumber \\
&\leq& \Big(\ell w e^{-p w / 2}\Big)^w \nonumber \\
&\leq& \Big(\ell w e^{-n^{\eps} /(4k')}\Big)^w \label{line:1} \\
&\leq& \Big(e^{-n^{\eps} /(5k')}\Big)^w \label{line:2} \\
&\leq& e^{-n^{1+\eps} /(10(k')^2)}. \label{line:3}
\end{eqnarray}
\noindent
In passing from \raf{line:1} to \raf{line:3}, the
reader may find it helpful to recall that $n= |V_i|$ (for $1\leq
i \leq k$), $k'=|V(C)|$, $\ell$ is between $2$ and $k$,
$w=\lceil n/(2k') \rceil$, and
$p=n^{\eps - 1}$, and that these estimates are valid for fixed $k'$
 and sufficiently large $n$. Of course, Lemma~\ref{lem:P_l}
follows if we take $\gamma= 1/(10(k')^2)$. \qed

\mn 
\textit{Proof 2 of Lemma~\ref{lem:P_l}}.  We use the Janson
Inequalities, (mainly) follow the notation of
\cite[Chapter~8]{Alon-Spencer2000}, and assume familiarity on the
readers' part. Here, $\Omega$ denotes the set of all potential
arcs (in $D_n$, as defined at the start of
Section~\ref{sect:proof2}) between the sets $U_{j_i}$, 
for $i=1,2...,\ell$, (introduced just prior to   our statement of
Lemma~\ref{lem:P_l}); each arc in $\Omega$ appears with
probability $p$.

Let $s$ be a (large) multiple of $\ell$; the value of $s$ will be
independent of $n$ and specified below. Now, let us enumerate
those cycles of $D_n$ that are of length $s$, and that cyclically
traverse $U_{j_1}, U_{j_2},..., U_{j_\ell}$ $s/ \ell$ times. For
$j\geq 1$, denote by $S_j$ the arc set of the $j$th such cycle and
by
$\EB_j$ the event that the arcs in $S_j$ all appear in $H$ (i.e.\
the cycle determined by $S_j$ is present in $H$). Let the random
variable $X$ count those $\EB_j$ that occur. Since $\Pr(X=0)$ (the
probability that there is no such cycle of length $s$) is an upper
bound for $P_{\ell}$ (the chance that $U(\ell)$ is acyclic), we
can bound $P_{\ell}$ by bounding $\Pr(X=0)$, and estimating the
latter quantity is exactly the purpose of Janson's Inequalities.
In the Janson paradigm, the value of $\Delta$ is defined by
\begin{equation}
\label{Janson-Delta-def1}
 \Delta \ddef \sum_{S_i \sim S_j}\Pr(\EB_i \cap \EB_j),
\end{equation}
where $S_i\sim S_j$ if the two cycles determined by $S_i$ and $S_j$
have at least one arc in common.

First, we find an upper bound for $\Delta$. Letting $i$ remain fixed, we
(rather crudely) obtain
\begin{equation}
\label{Delta-est0}
\Delta \leq n^s \sum_{j: S_i \sim S_j}\Pr(\EB_i \cap \EB_j),
\end{equation}
since each $|U_r|\leq n$ and each $|S_i|=s$. The sum on the right
side satisfies
\begin{equation}
\label{Delta-est1} \sum_{j: S_i \sim S_j}\Pr(\EB_i \cap \EB_j)
~~\leq~~ \sum_{r=1}^{s-1} \C{s}{r} p^{2s-r} w^{s-(r+1)}. %
\end{equation}
The binomial coefficient in (\ref{Delta-est1}) accounts for the
number of ways to choose the arcs of $S_i\cap S_j$, the power of
$p$ is $\Pr(\EB_j|\EB_i)\Pr(\EB_i)$, and finally, the power of
$w$ reflects the 
facts that each $U$-set has cardinality $w$ and, with $i$ fixed,
there are at most $s-(r+1)$ vertices in the $S_j$-cycle not
already in the $S_i$-cycle. %
Recalling that $w=\lceil n/(2k') \rceil$ (so that $w<n$), using
the gross bound $\C{s}{r}<2^s$, and replacing $p$ with
$n^{\eps-1}$, we find that
\[
\sum_{j: S_i \sim S_j}\Pr(\EB_i \cap \EB_j)~<~ 2^s
\sum_{r=1}^{s-1} p^{2s-r} n^{s-(r+1)}~=~ 2^s \sum_{r=1}^{s-1}
n^{2\eps s - s - r\eps -1}~<~ 2^s s  n^{2\eps s - s - \eps -1}.
\]
With (\ref{Delta-est0}), the last estimate yields
\begin{equation}
\label{Delta-est2} \Delta < 2^s s  n^{2\eps s - \eps -1}.
\end{equation}

Next, we find a lower bound for $\mu\ddef E[X]$. Since there are
$\ell$ $U$-sets, each containing $w$ vertices, and each ordered
choice of $s/\ell$ vertices from each (up to the choice of the
first vertex) contributes $1$ to $X$ with probability at least
$p^s$, we have
\[
\mu ~~\geq~~ \frac{1}{s}\C{w}{s/\ell}^\ell
\left[\left(\frac{s}{\ell}\right)!\right]^{\ell} p^s.
\]
Therefore,
\begin{equation}
\label{mu-est1} \mu ~\geq~ \frac{1}{s}\left(
\frac{w!}{(w-s/\ell)!} \right)^\ell p^s
    ~\geq~ \frac{1}{s}\left(w-\frac{s}{\ell} \right)^s p^s
    ~\geq~ \frac{1}{s}\left( \frac{n}{4k'} \right)^s n^{\eps s -s}
    ~=~ \frac{n^{\eps s}}{s(4k')^s}.
\end{equation}
\\
We distinguish two cases.
\\\\
\noindent
\textbf{Case 1:} $\Delta \geq \mu$.\\
Here, we have the hypotheses of the Extended Janson Inequality
(\cite[Theorem~8.1.2]{Alon-Spencer2000}), which, along with
our bounds (\ref{Delta-est2}), (\ref{mu-est1}) gives
\[
\Pr(X=0) ~\leq~ e^{-\mu^2/(2\Delta)} ~<~
e^{-n^{1+\eps}/(2s^3(32{k'}^2)^s)}.
\]

\noindent
\textbf{Case 2:} $\Delta < \mu$. \\
Now we have the hypotheses of the basic Janson Inequality
(\cite[Theorem~8.1.1]{Alon-Spencer2000}), which together with
(\ref{mu-est1}) gives
\[
\Pr(X=0) ~\leq~ e^{-\mu+\Delta/2} ~<~ e^{-\mu/2} ~\leq~
e^{-n^{\eps s}/(2s(4k')^s)}.
\]

\noindent
Let $s > 1+(1+\eps)/\eps$ be a multiple of $\ell$. Then the last
bound shows that
\[
\Pr(X=0) ~\leq~ e^{-n^{1+\eps}(n^{\eps}/(2s(4k')^s))} ~\leq~
e^{-n^{1+\eps}}.
\]

Since $s$ and $k'$ are constants (not depending on $n$), as is the
number $1$ (the coefficient of $n^{1+\epsilon}$ in the last expression), in
either case we see that
\[
P_{\ell} ~\leq~ \Pr(X=0) ~\leq~ e^{-\gamma n^{1+\eps}}
\]
for some constant $\gamma>0$. This gives us Lemma~\ref{lem:P_l}.
\qed

\bigskip
We return to our estimation of the expected number $N$ of bad
sequences in \raf{expected:N}, repeated here for convenience:
\[
N ~\leq~ \sum_{\ell=2}^{k} \C{k}{\ell} (\ell-1)!
{\C{n}{w}}^{\ell} P_{\ell}.
\]
Using Lemma~\ref{lem:P_l} to bound the factors $P_{\ell}$ in this
sum shows that for $n$ large enough,
\begin{equation}
\label{bound:N}
N ~\leq~ \sum_{\ell=2}^{k} \C{k}{\ell} (\ell - 1)!
{\C{n}{w}}^{\ell} e^{-\gamma n^{1+\eps}}
  ~<~ \sum_{\ell=2}^{k} \frac{e^{-n}}{2k}
 ~<~ \frac{e^{-n}}{2}.
\end{equation}

\noindent
From \raf{bound:N} and Markov's Inequality, we conclude that
\beq{badsequence}
\Pr(\exists\mbox{ a bad sequence}) < \frac{e^{-n}}{2}.
\enq

\noindent
Since $\phi$ fails to be an acyclic homomorphism
exactly when there exists a bad pair or there exists a bad sequence,
\raf{badpair} and \raf{badsequence} now show that
\[
\left|\DD\sm\Dii\right|\leq\left(\Pr(\exists\mbox{
bad pair})+\Pr(\exists\mbox{ bad
sequence})\right)\left|\DD\right| < e^{-n}\left|\DD\right|,
\]
which yields \raf{D2-lb}.
\qed

\section{Proof of Theorem~\protect{\ref{main-two}}}
\label{sect:proof3}

To obtain the conclusion of Theorem~\ref{main-two} (unique
$D$-colourability), we shall need to refine the deletion method
employed in the proof of Theorem~\ref{main-one}. We preserve the
earlier notation. Let $\Diii$ be the set of digraphs $H\in\Di$, in
which any two cycles of length less than $g$ are disjoint. Let
$\Diiii$ denote the set of those $H\in\DD$ with the property that
$H - A_1$ is uniquely $D$-colourable for any set $A_1$ of at most
$\nge$ independent arcs. (Here, we call a set $S \subseteq E(H)$
\DEF{independent} if no two arcs in $S$ have a vertex in common.)
Now we will show that \beq{D3-lb}
\left|\Diii\right|>\left(1-n^{-\eps/3}\right)\left|\DD\right| \enq
and \beq{D4-lb}
\left|\Diiii\right|>\left(1-e^{-n^\eps/6}\right)\left|\DD\right|.
\enq Since \raf{D3-lb} and \raf{D4-lb} imply that
$\Diii\cap\Diiii\neq\emp$ (for large enough $n$), we can choose a
digraph $H\in\Diii\cap\Diiii$. As $H\in\Diii\subseteq\Di$, we can
delete a set $A_1$ of at most $\nge$ independent arcs from $H$ so
that $D^*\ddef H - A_1$ has girth at least $g$, and $H\in\Diiii$
ensures that $D^*$ is uniquely $D$-colourable. Hence, to complete
the proof of Theorem~\ref{main-two}, it suffices to establish
\raf{D3-lb} and \raf{D4-lb}.

\bn \textit{Proof of\/} \raf{D3-lb}. For integers $\ell_1,\ell_2 <
g$, we follow \cite{Zhu-uniq} and call a digraph an
\DEF{$(\ell_1,\ell_2)$-double cycle} if it consists of a directed
cycle $C_{\ell_1}$ of length $\ell_1$ and a directed path of
length $\ell_2$ joining two (not necessarily distinct) vertices of
$C_{\ell_1}$; such a digraph contains $\ell_1+\ell_2-1$ vertices
and $\ell_1+\ell_2$ arcs. Let $\DC$ denote the set of digraphs in
$\DD$ containing an $(\ell_1,\ell_2)$-double cycle for some
$\ell_1,\ell_2 < g$. Notice that $\Di\sm\Diii\subseteq\DC$, whence
\beq{Dprime-est} 
\left|\Di\sm\Diii\right|\leq\left|\DC\right|,
\enq 
so we can obtain a lower estimate for $\left|\Diii\right|$ by
estimating $|\DC|$.

For fixed $\ell_1,\ell_2 < g$, the expected number
$N(\ell_1,\ell_2)$ of $(\ell_1,\ell_2)$-double cycles in a
digraph $H\in\DD$ is less than
\[
\ell_1(kn)^{\ell_1}(kn)^{\ell_2-1} p^{\ell_1 + \ell_2},
\]
since there are fewer than $\ell_1(kn)^{\ell_1}(kn)^{\ell_2-1}$
ways of choosing such a double cycle $Y$ with $V(Y)\subseteq V$,
and each such $Y$ exists with probability $0$ or
$p^{\ell_1+\ell_2}$. Since $p=n^{\eps-1}$ we have
\[
N(\ell_1,\ell_2) <
\ell_{1}k^{\ell_1+\ell_2}n^{\eps(\ell_1+\ell_2)}n^{-1}.
\]
Since $\eps(\ell_1 + \ell_2) \leq 2g\eps<1/2$, for large enough
$n$ we have
\[
\sum_{\mystack{2\leq\ell_1<g}{1\leq\ell_2<g}}N(\ell_1,\ell_2)<n^{-1/2}.
\]
Markov's Inequality now shows that
\[
\left|\DC\right|<n^{-1/2}\left|\DD\right|,
\]
so from \raf{Dprime-est} we obtain
\[
\left|\Diii\right| > \left|\Di\right|-n^{-1/2}\left|\DD\right|,
\]
and \raf{D1-lb} gives \raf{D3-lb}. \qed

\medskip
\noindent 
\textit{Proof of\/} \raf{D4-lb}. We will argue that
$|\DD\sm\Diiii|<e^{-n^{\eps}/6}|\DD|$. If $H\in\DD\sm\Diiii$, then
there is a set $A_1$ of at most $\nge$ independent arcs of $H$ so
that $H-A_1$ admits an acyclic homomorphism $h$ to $D$ that is
not the composition $\sigma\circ c$ of the natural homomorphism
$c\colon H-A_1\to D$ (sending $V_i$ to $i$) with an automorphism
$\sigma$ of $D$. As in the proof of \raf{D2-lb}, we can define a
function $\phi\colon V(D)\to V(D)$ such that $\left|V_i\cap
h^{-1}(\phi(i))\right|\geq n/k$ for each $i\in V(D)$.  %
Let us first suppose that $\phi$ is not an automorphism of $D$. By
hypothesis, $D$ is a core, so any acyclic homomorphism of $D$ to
itself must be an automorphism. It follows that $\phi$ is not an
acyclic homomorphism. Therefore, there is an arc $ij\in E(D)$ such
that $\phi(i)\phi(j)\not\in E(D)$, or there is a vertex $i\in
V(D)$ such that $\phi^{-1}(i)$ is not acyclic. Notice that the
current arrangement is analogous to the one in the second
paragraph in the proof of \raf{D2-lb}. Repeating the earlier
argument, with $D$ in the place of $C$ and $k$ in the role of
$k'$, we find that most $H\in\DD$ do not fall into the present
case. More precisely, we reach the following conclusion:

\mn
\textit{At least $(1-e^{-n})|\DD|$ digraphs $H\in\DD$ have the
  property that for any set $A_1$ of at most $\nge$ arcs
  (independent or otherwise), the digraph $H-A_1$ cannot be
  $D$-coloured so that $\phi$ is not an automorphism of $D$.
}

\mn Thus, in this case,
$|\DD\sm\Diiii|<e^{-n}|\DD|<e^{-n^{\eps}/6}|\DD|$, and \raf{D4-lb}
is proved.

From now on, we treat the case when $\phi$ is an automorphism of
$D$. Without loss of generality, we may assume that $\phi$ is the
identity, i.e., that \beq{preimage} \left|V_i\cap
h^{-1}(i)\right|\geq n/k \mbox{ for each } i\in V(D). \enq We may
assume further that \beq{upperbound} \left|V_j\cap
h^{-1}(i)\right| < n/k \mbox{ for all }j\neq i. \enq (Otherwise,
we can redefine $\phi(i)$ to be equal to $j$ and fall into the
case where $\phi$ is not an automorphism.) %

Since $h$ is not the composition $\sigma\circ c$ of the natural
homomorphism $c\colon H-A_1\to D$ (sending $V_i$ to $i$) with an
automorphism $\sigma$ of $D$, there must be a pair $\{i,j\}$ of
distinct vertices of $D$ such that $ V_{j} \cap h^{-1}(i)\neq
\emp$. Let $\{i_0,j_0\}$ be such a pair that maximizes
$|V_{j_0}\cap h^{-1}(i_0)|$. Consider the map $\phi'\colon V(D)
\to V(D)$ such that
\[
\phi'(x) \ddef \left\{
\begin{array}{ll}
x~(=\phi(x)) &\mbox{if $x\neq j_0$} \\
i_0 & \mbox{if $x=j_0$}.
\end{array}\right.
\]
Clearly $\phi'$ is not a bijection, and since $D$ is a core, it
cannot be an acyclic homomorphism. There are two possibilities.

\mn \emph{Case} 1: Both $j_0 i_0$ and $i_0 j_0$ are arcs of $D$
(so $\phi'^{-1}(i_0)$ is not acyclic).

\mn \emph{Case} 2: There exists $v\in V(D)$ such that $vj_0$ is an
arc of $D$ but $v i_0 $ is not, or $j_0 v$ is an arc of $D$ but
$i_0 v$ is not.

\bn We will show that in either case,
$|\DD\sm\Diiii|<e^{-n^{\eps}/6}|\DD|$.

\mn \textbf{Case 1:} Our choice of $\{i_0, j_0\}$ ensures that
$h^{-1}(i_0)\cap V_{j_0}\neq\emp$. Let $x\in h^{-1}(i_0)\cap
V_{j_0}$, and consider the (nonrandom) subdigraph $\widehat{D_n}$
of $D_n$ induced by $\{x\}\cup(h^{-1}(i_0)\cap V_{i_0})$. As
$V_{i_0}$ induces no cycles, all cycles of $\widehat{D_n}$ must
include $x$, and since the arcs of $A_1$ are independent, at most
one such arc is incident with $x$. Furthermore, the constraint on
the size of $A_1$ and our choice of $\eps$ (smaller than $1/4g$)
give
\[
|A_1| \leq \lceil n^{g\eps}\rceil < \lceil n^{1/4}\rceil \ll \frac{n}{k}.
\]
Because $|h^{-1}(i_0)\cap V_{i_0}|\geq n/k$ (cf.\ \raf{preimage}),
there must be a subset $U \sub h^{-1}(i_0)\cap V_{i_0}$ of
cardinality $\lfloor n/2k \rfloor$ such that the (random)
subdigraph induced by $\{x\}\cup \,U$ contains no arcs of $A_1$
and moreover is acyclic (since $h^{-1}(i_0)$ is acyclic). To show
that this is unlikely, we first estimate the expected number $M$
of ways to select a vertex $x\in V_{j_0}$ and a subset $U \sub
V_{i_0}$ of cardinality $\lfloor n/2k \rfloor$ so that the
subdigraph $H_{x,U}$ of $H$ that they induce is acyclic and no arc
of $A_1$ is incident with a vertex in $U$. If $P_{x,U}$ denotes
the probability that $H_{x,U}$ is acyclic, then \beq{eq:T} M \leq
n \C{n}{\lfloor n/2k \rfloor}P_{x,U} < n^{n} P_{x,U}. \enq

In order to estimate $P_{x,U}$, we again employ the Janson
Inequalities (cf.\ \cite[Chapter 8]{Alon-Spencer2000}). Now
$\Omega$ denotes the set of all potential arcs in the subdigraph
$D_{x,U}'$ of $D_n$ induced by $\{x\} \cup U$; each arc in
$\Omega$ appears in $H_{x,U}$ with probability $p$. Let
$\ell>(2+\eps)/\eps$ be a fixed integer. Let us index those cycles
of $D_{x,U}'$ (with the positive integers) that are of length
$\ell+1$ in $D_{x,U}'$. For $j\geq 1$, let $S_j$ be the arc-set of
the $j$th such cycle and $\EB_j$ be the event that the arcs in
$S_j$ all appear (i.e.\ the cycle determined by $S_j$ is present
in $H_{x,U}$). Let $X$ count the $\EB_j$ that occur; since
$\Pr(X=0)$ is an upper bound for $P_{x,U}$, we can bound $P_{x,U}$
by bounding $\Pr(X=0)$.

As in \raf{Janson-Delta-def1}, Janson's $\Delta$ is given by
\[
 \Delta \ddef \sum_{S_i \sim S_j}\Pr(\EB_i \cap \EB_j).
\]

\noindent
Since there are at most
$\C{\left\lfloor n/2k\right\rfloor}{\ell}<n^{\ell}$ cycles
$S_j$, if we fix an $S_i$ to maximize
$\sum_{j: S_j \sim S_i}\Pr(\EB_i\cap \EB_j)$, then
\begin{equation}
\label{new-Delta-est0}
\Delta \leq n^{\ell} \sum_{j: S_j \sim S_i} \Pr(\EB_i \cap\EB_j).
\end{equation}
Now we sum over the number $r$ of common arcs an $S_j$ can
have with $S_i$;
this fixes at least $r+1$ vertices of $S_j$. Thus,
\[
\sum_{j: S_j \sim S_i} \Pr(\EB_i \cap \EB_j)
~\leq~ \sum_{r=1}^{\ell}
\C{\ell+1}{r}{\left\lfloor\frac{n}{2k}\right\rfloor}^{\ell-r-1} p^{2(\ell+1)-r}.
\]
Using the crude upper estimates $\C{\ell+1}{r}<2^{\ell+1}$   and
$\lfloor n/2k \rfloor<n$, and replacing $p$ with $n^{\eps-1}$, we
obtain
\[
\sum_{j: S_j \sim S_i}\Pr(\EB_i\cap\EB_j) ~<~ 2^{\ell+1}
\sum_{r=1}^{\ell} (np)^{\ell-r-1}p^{\ell+3} ~<~ 2^{\ell+1} \ell
(np)^{\ell-2}p^{\ell+3} ~=~ 2^{\ell+1} \ell
n^{2\eps\ell+\eps-\ell-3}.
\]

\noindent
This and \raf{new-Delta-est0} now give
\begin{equation}
\label{new-Delta-est1} \Delta \leq 2^{\ell+1} \ell
n^{2\eps\ell+\eps-3}.
\end{equation}

We also need to find a lower bound for $\mu\ddef E[X]$. Since the
arcs of $D_{x,U}'$ within $U$ are acyclically oriented, each
choice of $\ell$ vertices within $U$ determines exactly one
potential $(\ell+1)$-cycle (viz., through $x$). It follows that
\begin{equation}
\label{new-mu-est0}
\mu ~=~ \C{\lfloor n/2k \rfloor}{\ell}p^{\ell+1}
~>~
\left(\frac{\lfloor n/2k\rfloor}{\ell}\right)^{\ell} p^{\ell+1}
~>~ \frac{n^{\eps \ell + \eps - 1}}{(4k\ell)^{\ell}}.
\end{equation}

\bn As in the proof of Theorem~\ref{main-one}, we have two
subcases.

\noindent
\textbf{Subcase 1(i):} $\Delta \geq \mu$.\\
Again, we have the hypotheses of the Extended Janson Inequality
(\cite[Theorem~8.1.2]{Alon-Spencer2000}), which, along with
(\ref{new-Delta-est1}) and (\ref{new-mu-est0}) gives
\[
\Pr(X=0)~\leq~ e^{-\mu^2/(2\Delta)} ~<~ e^{-n^{1+\eps}/(\ell2^{\ell+2}(4k\ell)^{2\ell})}
\deff e^{-\beta n^{1+\eps}},
\]
where $\beta$ is the (positive) constant (not depending on
$n$) absorbing the denominator in the preceding exponent.

\noindent
\textbf{Subcase 1(ii):} $\Delta < \mu$. \\
Here, we have the hypotheses of the Janson Inequality
(\cite[Theorem~8.1.1]{Alon-Spencer2000}), which, with the help
of (\ref{new-mu-est0}) gives
\[
\Pr(X=0) ~\leq~ e^{-\mu+\Delta/2} ~<~ e^{-\mu/2} ~<~
e^{-n^{\eps\ell+\eps-1}/(2(4k\ell)^{\ell})}.
\]

\noindent
Recalling our choice of $\ell>(2+\eps)/\eps$, we see that
\[
\Pr(X=0)  ~<~  e^{-n^{1+2\eps}/(2(4k\ell)^{\ell})} ~<~ e^{-n^{1+\eps}}.
\]

In either subcase, we have that $P_{x,U}\leq \Pr(X=0)< e^{-\beta
n^{1+\eps}}$ (since $\beta<1$), and returning to \raf{eq:T}, we
have
\[
M ~<~ n^{n} P_{x,U} ~<~ n^ne^{-\beta n^{1+\eps}}
  ~=~ \left(ne^{-\beta n^\eps}\right)^n
  ~<~ e^{-\beta n^{1+\eps}/2}.
\]

\noindent 
By Markov's Inequality, the probability that there
exists such an $\{x\}\cup U$ (that induces an acyclic subdigraph)
is less than $e^{-\beta n^{1+\eps}/2}<e^{-n^{\eps}/6}$, and so in
Case~1, $|\DD\sm\Diiii|<e^{-n^{\eps}/6}|\DD|$, as desired.

\mn \textbf{Case 2:} By the hypothesis of this case, there is a
vertex $v$ such that either $vj_0\in E(D)$ and $vi_0\not\in E(D)$,
or $j_0v\in E(D)$ and $i_0v\not\in E(D)$. We will consider the
first of these; the second one yields to similar reasoning. Let us
recall that we chose a pair $\{i_0,j_0\}$ of distinct vertices of
$D$ so as to maximize $b:=|V_{j_0}\cap h^{-1}(i_0)|\neq 0$.

\mn \textbf{Claim:} Every vertex $z\in V(D)\sm\{i_0\}$ satisfies
$\left|V_z \cap h^{-1}(z)\right| \geq n-(k-1)b$.

\mn \textit{Proof of claim.} Otherwise, some $z\neq i_0$ satisfies
$\left|V_z \cap h^{-1}(z)\right|<n-(k-1)b$. By the pigeonhole
principle, there is some $u\neq z$ such that $\left|V_z \cap
h^{-1}(u)\right| > b$, but this contradicts our choice of
$\{i_0,j_0\}$. \qed

Using the claim, we see that there are sets $U_v\subseteq V_v\cap
h^{-1}(v)$ and $U_{j_0}= V_{j_0}\cap h^{-1}(i_0)$ with
$|U_v|=n-(k-1)b$ and $|U_{j_0}|=b$. Since $h\colon H - A_1\to D$
is an acyclic homomorphism and $vi_0\not\in E(D)$, there are at
most $\min\{b,\nge\}$ independent arcs from a vertex in $U_v$ to
one in $U_{j_0}$. We now estimate the expected number $L(b)$ of
pairs $U_v'\subseteq V_v$, $U_{j_0}'\subseteq V_{j_0}$ with
$|U_v'|=n-(k-1)b=n-(k-1)|U_{j_0}'|$, and at most $\min\{b,\nge\}$
arcs from $U_v'$ to $U_{j_0}'$.

For $b<n/k$ (cf.\ $\raf{upperbound}$) and $s\leq\min\{b,\nge\}$,
denote by $L(b,s)$ the expected number of pairs $U_v'\subseteq
V_v$, $U_{j_0}'\subseteq V_{j_0}$,
$|U_v'|=n-(k-1)b=n-(k-1)|U_{j_0}'|$, and exactly $s$ arcs joining
a vertex in $U_v'$ to one in $U_{j_0}'$. Then
\begin{eqnarray*}
L(b,s)&<&\C{n}{n-(k-1)b}\C{n}{b}\C{(n-(k-1)b)b}{s}p^s(1-p)^{(n-(k-1)b)b-s}\\
&<&n^{(k-1)b}n^b (nb)^s n^{s(\eps -1)} e^{-bn^\eps + n^{\eps-1}((k-1)b^2+s)}\\
&<& b^s n^{\eps s} n^{kb} e^{-(bn^\eps)/2}\\
&=& b^s n^{\eps s} (n^{k} e^{-n^\eps/2})^b\\
&<& b^s n^{\eps s} e^{-(bn^\eps)/3}\\
&<&e^{-n^\eps / 4}.
\end{eqnarray*}
\noindent
Letting
$L(b)=\sum_{s\leq\min\{b,\nge\}}L(b,s)<\nge e^{-n^\eps/4}<e^{-n^\eps/5}$,
we find that
\[
\sum_{1\leq b < n/k}L(b) < (n/k) e^{-n^\eps/5} < e^{-n^\eps/6}.
\]

This completes the discussion for the case when  $vj_0\in E(D)$
and $vi_0\not\in E(D)$; an identical argument gives the same upper
bound in the case when $j_0v\in E(D)$ and $i_0v\not\in E(D)$. 
Thus in Case~2, we also arrive at $|\DD\sm\Diiii|<e^{-n^{\eps}/6}|\DD|$.

\medskip
Combining the estimates obtained above and applying Markov's
Inequality finally yields \raf{D4-lb} and therefore completes the
proof of Theorem~\ref{main-two}. \qed

\section{The circular chromatic number}
\label{sect:chic}

We turn now to the implications of Theorem~~\ref{main-two} for 
circular colouring digraphs. The concept of the digraph circular
chromatic number $\chic$, defined below, generalizes the circular
chromatic number for undirected graphs. The theory of
the graph invariant, as of 2001, was surveyed in
\cite{Zhu-surv}. The digraph version was introduced in
\cite{BFJKM2002}, where it was proved, via Lemma~\ref{lem:chic}
below, that $\chic$ assumes all rational values at least one.
(Note that the digraphs of Lemma~\ref{lem:chic} do not generally have 
large girth.) The same article also established the following analogue 
of the Erd\H{o}s' theorem introducing the present paper: 
there exist digraphs with arbitrarily large girth 
and arbitrarily large circular chromatic number 
(this is the result to which we alluded
immediately following the statement of Theorem~~\ref{main-one}).
Our main result here, Theorem~\ref{thm:colouring}, provides 
a common generalization and strengthening of these basic results.
It shows that the `all conceivable rationals' property of $\chic$
holds even for digraphs of arbitrarily large girth and even
demanding a certain uniqueness of the colouring. 

Let $d\ge1$ and $k\ge d$ be integers. Let $C(k,d)$ be the digraph with
vertex set $\ZZ_k = \{0,1,\dots,k-1\}$ and arcs
\[
    E(C(k,d)) = \{ ij\mid j-i\in \{d,d+1,\dots,k-1\} \},
\]
where the subtraction is considered in the cyclic group $\ZZ_k$ of integers
modulo $k$.

Acyclic homomorphisms into $C(k,d)$ are an important concept
because of their relation to the circular chromatic number of
digraphs; cf.~\cite{BFJKM2002}. An acyclic homomorphism of a
digraph $D$ into $C(k,d)$ is called a \DEF{$(k,d)$-colouring} of
$D$. It is shown in \cite{BFJKM2002,Mo} that there is a rational
number $q\in \QQ$ such that $D$ has a $(k,d)$-colouring if and only
if $k/d \ge q$. This value $q$ is denoted by $\chic(D)$ and
called the \DEF{circular chromatic number} of $D$. For 
$q \in \QQ^+$, let $S_q$ denote the circle of perimeter $q$
(centred, say, at the origin of $\RR^2$).
We define a \DEF{circular $q$-colouring} of $D$ to be a map
$\phi:V(D) \to S_q$ such that for every $xy \in E(D)$, 
with $\phi(x) \neq \phi(y)$, the
distance $d_S(\phi(x), \phi(y))$ from $\phi(x)$ to $\phi(y)$ in
the clockwise direction around $S_q$ is at least $1$, and for every
$p \in S_q$, the preimage $\phi^{-1}(p)$ induces an acyclic
subdigraph of $D$. If $\phi$ is a
circular $q$-colouring, we say that the arc $xy \in E(D)$ is
\DEF{tight} whenever $d_S(\phi(x),\phi(y)) \leq 1$ (in which case this
distance is either 1 or 0). A cycle in $D$ consisting of tight arcs is
called a \DEF{tight cycle} for the circular $q$-colouring $\phi$.
Note that every tight cycle contains an arc $xy$ such that
$d_S(\phi(x),\phi(y))=1$. We will use the following results,
respectively from \cite{Mo} and \cite{BFJKM2002}.

\begin{lemma} 
\label{lem:circol}
If $\chi_c(D) = q$, then every circular $q$-colouring of $D$ has a
tight cycle.
\end{lemma}

\begin{lemma}
\label{lem:chic} $\chic(C(k,d))=k/d$.
\end{lemma}

Lemmas~\ref{lem:circol} and \ref{lem:chic} imply the following fact.

\begin{prop}
\label{prop:core} If $k$ and $d$ are integers with $1\le d \leq k$,
then $C(k,d)$ is a core if and only if $k$ and $d$ are relatively
prime.
\end{prop}

\proof Let $C=C(k,d)$ and $V=V(C)$. If $r\ddef\gcd(k,d)>1$, then the
mapping $\phi\colon V\to V$ given by $\phi(i) \ddef r\lfloor
i/r\rfloor$ is easily seen to be an acyclic homomorphism $C \to C$
that is not surjective. By Lemma~\ref{lem:core}, $C(k,d)$ is not a core.

For the converse, assume that $k$ and $d$ are relatively prime, and
let $\phi\colon V\to V$ be an acyclic homomorphism. Define a map
$\varphi\colon C(k,d) \to S_{k/d}$ as follows. Let 
$s_0, s_1,\ldots,s_{k-1}$ be points on $S_{k/d}$ that appear on the circle
consecutively at distance $1/d$ apart. For $0\leq i\leq k-1$, we 
set $\varphi(i) \ddef s_{\phi(i)}$. Since $\phi$ is an acyclic
homomorphism, it is easily verified that $\varphi$ is a
circular $\frac{k}{d}$-colouring of $C(k,d)$. 
By Lemmas~\ref{lem:circol} and \ref{lem:chic},
$\varphi$ has a tight cycle $C_0 = v_1v_2\cdots v_{\ell}v_1$ in 
$C(k,d)$. We
may assume that $\varphi(v_1) = s_0$. The images $\varphi(v_1),
\varphi(v_2),\ldots,\varphi(v_{\ell}),\varphi(v_1)$ must take 
consecutive values $s_0, s_d, s_{2d}, s_{3d},\ldots$ 
(each possibly repeated several
times), with the indices taken modulo $k$, and end up at $s_0$.
Since $k$ and $d$ are relatively prime, the sequence $s_0, s_d,
s_{2d},\ldots$ must exhaust all the elements in the set $\{s_0, s_1,
\ldots, s_{k-1} \}$. This shows that $\phi$ is surjective; whence, by
Lemma~\ref{lem:core}, $C(k,d)$ is a core.  \qed 

Proposition~\ref{prop:core} and Theorem~\ref{main-two} together
yield an immediate consequence, Corollary~\ref{cor:colouring}, 
that we now state in a slightly expanded (and more precise) 
form:

\begin{thm}
\label{thm:colouring} If $k$ and $d$ are relatively prime
integers with $1 \leq d \leq k$, then for every positive integer $g$,
there exists a uniquely $C(k,d)$-colourable digraph of girth at
least $g$ (and with circular chromatic number equal to $k/d$).
\end{thm}

\nin
The last claim of Theorem~\ref{thm:colouring} follows from
the next result, an
analogue of \cite[Theorem 3]{Zhu-uniq}.

\begin{thm}
If $D$ is a uniquely $C(k,d)$-colourable digraph, 
then $\chic(D)=k/d$.
\end{thm}

\proof 
Since $D$ is $C(k,d)$-colourable, we have $\chic(D) \leq k/d$.
Suppose, for a contradiction, that
$\chic(D) = k'/d' < k/d$. Define $d^{*}\ddef dd'$, $m \ddef kd'$ and
$m'\ddef k'd$ so that 
$m'/d^{*}=k'/d'<k/d=m/d^{*}$.
Now, let $\phi'$ be an $(m', d^{*})$-colouring of $D$. 
Using the idea in the proof of Proposition~\ref{prop:core}, we can
define a circular $\frac{m'}{d^{*}}$-colouring $\varphi$ of 
$C(m',d^{*})$ so that $\varphi\circ\phi'$ is such a colouring of $D$.
Since $\chic(D)=k'/d'=m'/d^{*}$, Lemma~\ref{lem:circol} implies 
that $\varphi\circ\phi'$ has a tight cycle in $D$. Choosing a tight
arc $xy \in E(D)$ for $\varphi\circ\phi'$ yields an arc $xy$ of $D$
such that $\phi'(x)$ and $\phi'(y)$ are
separated by $d^{*}$ units in the clockwise direction around
$C(m',d^{*})$.
Without loss of generality, we
may assume that $\phi'(y) = 0$ and $\phi'(x) = m'-d^{*}$. We
define an $(m, d^{*})$-colouring $\psi$ as follows: $\psi(v) \ddef
\phi'(v)$ if $\phi'(v) < m' - d^{*}$ and $\psi(v) \ddef \phi'(v) + m -
m'$ otherwise. It is easily verified that $\psi$ is indeed an 
$(m,d^{*})$-colouring of $D$. Next, we define 
$\bar{\psi} \colon V(D) \to V(C(k,d))$ by 
$\bar{\psi}(v) \ddef \left\lfloor \psi(v)/d' \right\rfloor$. 
(In this formula---and hereafter---we view the
vertices $\psi(v)$ of $C(m,d^{*})$ as integers between $0$ and
$m-1$.) 
Since $\lfloor\cdot/d'\rfloor\colon V(C(m,d^{*}))\to V(C(k,d))$
defines an acyclic homomorphism, and such maps compose
(cf.~\cite{BFJKM2002}), it is not hard to check that $\bar{\psi}$ is a
$(k,d)$-colouring of $D$. Similarly, we define 
$\bar{\phi} \colon V(D) \to V(C(k,d))$ by 
$\bar{\phi}(v) \ddef \left\lfloor \phi'(v)/d'\right\rfloor$. 
As in the case of $\bar{\psi}$, it is easy to check 
that $\bar{\phi}$ is a $(k,d)$-colouring of $D$. 
We claim that $\bar{\phi}$ and $\bar{\psi}$ do
not differ by an automorphism of $C(k,d)$. Note that $\bar{\phi}(y) =
\bar{\psi}(y) = 0$; therefore, it suffices to show that
$\bar{\phi}(x) \neq \bar{\psi}(x)$. Now, 
\[ 
\bar{\phi}(x) = \left\lfloor \frac{m'-d^{*}}{d'} \right\rfloor 
= \left\lfloor \frac{d(k'-d')}{d'} \right \rfloor \qquad \mbox{while} \qquad
\bar{\psi}(x)= \left\lfloor \frac{m-d^{*}}{d'} \right\rfloor = k-d.
\]
Since  $k'/d'< k/d$, we have 
$d(k'-d')/d' < k-d$, and it follows that $\bar{\phi}(x) <
\bar{\psi}(x)$. This implies that $\bar{\phi}$ and $\bar{\psi}$
are $(k,d)$-colourings of $D$ that do not differ by an 
automorphism of $C(k,d)$. Hence, $D$ is not uniquely 
$C(k,d)$-colourable, a contradiction.
\qed

\end{document}